\documentclass[11pt,a4paper]{amsproc}

\newcommand{\editorsnote}{\relax}

\newtheorem{theorem}{Theorem}[section]

\newtheorem{lemma}{Lemma}[section]
\newtheorem{proposition}{Proposition}[section]
\newtheorem{definition}{Definition}[section]

\newcommand{\kpwcite}[1]{\cite{#1}}

\newtheorem{problem}{Problem}[section]

\newcommand{\range}{\operatorname{range}}

\newcommand{\Gr}{\operatorname{Gr}}

\newcommand{\cyl}{\textup{cyl}}

\newcommand{\detz}{\operatorname{det_\zeta}}
\newcommand{\detf}{\operatorname{det_F}}
\newcommand{\KPW}{Krzysztof P. Wojciechowski}
\newcommand{\bvp}{boundary value problem}

\newcommand{\sGr}{\operatorname{Gr}_\infty^*}

\newcommand{\etab}{\widetilde\eta}

\newcommand{\cinfz}[1]{C_0^\infty(#1)}
\newcommand{\cinf}[1]{C^\infty(#1)}
\newcommand{\ginf}[1]{\Gamma^\infty(#1)}

\newcommand\Calderon{Calder{\'o}n}
\newcommand{\End}{\operatorname{End}}
\newcommand{\Hom}{\operatorname{Hom}}
\newcommand{\Aut}{\operatorname{Aut}}

\newcommand{\dvol}{\operatorname{dvol}}
\newcommand{\dom}{\operatorname{dom}}

\newcommand{\mlmod}{\operatorname{mod}}

\newcommand{\bigsetdef}[2]{\bigl\{ #1 \,\bigm|\, #2\bigr\}}

\newcommand{\scalar}[2]{\langle #1,#2\rangle}
\newcommand{\scalarround}[2]{\left(  #1,#2 \right)}

\newcommand\ga{\alpha}

\newcommand\pl{\partial}

\newcommand\gl{\lambda}

\newcommand\eps{\varepsilon}

\newcommand{\C}{\mathbb{C}}
\newcommand{\R}{\mathbb{R}}

\newcommand{\Z}{\mathbb{Z}}

\newcommand\cE{\mathcal{E}}

\newcommand\cP{\mathcal{P}}

\newcommand\coker{\operatorname{coker}}

\newcommand\ind{\operatorname{ind}}

\renewcommand\Re{\operatorname{Re}}

\newcommand\sign{\operatorname{sign}}
\newcommand\spec{\operatorname{spec}}

\newcommand\tr{\operatorname{tr}}

\begin{document}

\title[The mathematical work of Krzysztof P. Wojciechowski]{Selected Aspects of the Mathematical Work of Krzysztof P. Wojciechowski}

\author{Matthias Lesch}

\thanks{{S}upported partially
by \uppercase{S}onderforschungsbereich/\uppercase{T}ransregio 12 
``\uppercase{S}ymmetries and \uppercase{U}niversality in
\uppercase{M}esoscopic \uppercase{S}ystems" 
(\uppercase{B}ochum--\uppercase{D}uisburg/\uppercase{E}ssen--\uppercase{K}\"oln--\uppercase{W}arszawa)}

\address {Mathematisches Institut\\
Beringstr. 1\\
53115 Bonn, Germany
}
\email{lesch@math.uni-bonn.de}

\dedicatory{Dedicated to Krzysztof P. Wojciechowski on his 50th birthday}


\begin{abstract}
To honor and to please our friend Krzysztof P. Wojciechowski I will review the
milestones of his mathematical work. This will at the same time be a tour
of Analysis and Geometry of Boundary Value Problems. Starting in
the 80s I will discuss the spectral flow and the general linear conjugation problem,
the \Calderon\ projector and the topology of space of elliptic boundary problems.
The theme of the 90s is the eta invariant. The paper with Douglas
was fundamental for establishing spectral invariants 
for manifolds {\em with} boundary and for 
the investigation of the behavior of spectral invariants
under analytic surgery. This was so influential that many different proofs of
the gluing formula for the eta-invariant were published.
Finally turning to the new millennium we will look at the zeta--determinant. Compared to
eta this is a much more rigid spectral invariant which is technically challenging.
\end{abstract}

\maketitle


\section{Introduction}
\subsection{The framework and the problem}

To begin with let us describe in general terms the problems to which \KPW\ has
contributed so much in the last 25 years. 

Let $X$ be a compact smooth Riemannian manifold with boundary $\Sigma=\partial
X$.
Furthermore, let $E,F$ be hermitian vector bundles over $X$ and let 
\begin{equation}
    D:\ginf{X,E}\longrightarrow \ginf{X,F}
\end{equation}
be an elliptic differential operator: $\ginf{X,E}$ denotes the spaces
of smooth sections of the bundle $E$.

In this situation some natural questions occur:

\par\bigskip
1. \textit{What are appropriate boundary conditions for $D$ on $X$?}
\par\medskip

This question is absolutely fundamental since without imposing boundary
conditions we cannot expect $D$ to have any reasonable spectral theory.

A boundary condition is given by a pseudo--differential operator
\begin{equation} 
    P:\ginf{\Sigma,E}\longrightarrow \ginf{\Sigma,E}  
\end{equation}
of order $0$.\footnote{One could think of more general definitely nonlocal boundary operators, but
in this paper we will content ourselves to pseudo--differential boundary conditions.}
The \emph{realization} $D_P$ of the boundary condition given by $P$ is the
differential expression $D$ acting on the domain 
\begin{equation}\label{domain}
    \dom(D_P):=\bigsetdef{u\in L_1^2(X,E)}{P(u|\Sigma)=0}.
\end{equation}

Since $D$ is elliptic what one should expect naturally for $P$ to be
"appropriate"
is that \emph{elliptic regularity} holds. That is if $Du\in
L^2_s(X,E)$
     \footnote{We denote the space of sections of $E$ which are of
                 Sobolev order $s$ by $L^2_s(X,E)$.}
is of Sobolev order $s\ge 0$ and if $P(u|\Sigma)=0$ then $u\in L^2_{s+d}(X,E)$
is already of Sobolev order $s+d$, where $d$ denotes the order of $D$.

\par\bigskip
2. \textit{What is the structure of the space of all (nice) boundary
conditions and how do spectral invariants of $D_P$ depend on the boundary condition?}
\par\medskip

These problems are the Leitfaden of \KPW's work.
If we are given a realization $D_P$ of a nice boundary value problem we can do
spectral theory and study the basic spectral invariants of $D_P$. 
We will see that the question in the headline
leads to interesting and delicate analytical
problems. Let us specify the kind of spectral invariants we mean here.

The most basic spectral invariant of the Fredholm operator $D_P$ is its index
\begin{equation} 
    \ind D_P=\dim \ker D_P-\dim \coker D_P.
\end{equation}

More rigid (and analytically more demanding) spectral invariants are derived
from the heat trace
\begin{equation} 
   \tr \bigl(e^{-tD_P^2}\bigr)=\sum_{\lambda\in\spec D_P\setminus\{0\}} e^{-t\lambda^2},
\end{equation}
where $D_P$ is now assumed to be self--adjoint, via Mellin transforms. 
The most important examples are
the $\eta$--invariant
\begin{equation}\label{eta}
   \eta(D_P)=\left[\frac{1}{\Gamma(\frac{s+1}{2})}\int_0^\infty
	  t^{(s-1)/2}\tr\bigl(D_Pe^{-tD_P^2}\bigr)dt\right]_{s=0}
\end{equation}
and the $\zeta$--determinant
\begin{equation}\label{zeta}
      \log\detz(D)=-\frac{d}{ds}\left[\frac{1}{\Gamma(s)}\int_0^\infty
	t^{s-1}\tr \bigl(e^{-tD_P^2}\bigr)dt \right]_{s=0}.
\end{equation}

The existence of these invariants is highly non--trivial
since it depends on the meromorphic continuation of the 
right hand side of \eqref{eta} and \eqref{zeta}.

In the following sense the index is the least rigid and the 
$\zeta$--determinant is the most rigid of these three
invariants. In order not to get into too much technicalities
assume for the moment that $D(s)_{a\le s\le b}$ is a smoothly varying family
of elliptic operators on a closed manifold.\footnote{Here smoothly varying
means that all coefficients depend smoothly on the parameter.}

The index is insensitive to small perturbations of the operator.
Hence $\ind D(s)$ will not depend on $s$ at all. The variation
of the $\eta$--invariant is easy to understand. First of all
the \emph{reduced} $\eta$--invariant 
\begin{equation}
\widetilde \eta(D(s))=\frac 12(\dim\ker D(s)+\eta(D(s))
\end{equation}
has only integer jumps and the total number
of jumps equals the \emph{spectral flow} of the family $D(s)$ over
the interval [a,b]. The variation of $\widetilde\eta(D(s))\mlmod \Z$
is \emph{local} in the sense that $\frac{d}{ds}(\widetilde\eta(D(s)\mlmod \Z)$
is the integral of a density which is a local expression in terms of the
coefficients of the operator and its derivatives, cf. Gilkey \kpwcite{Gil:ITH2E}, Sec. 1.13.
The variation of the $\zeta$--determinant is more complicated and
depends on global data.

\begin{sloppy}
It is therefore most natural that the early work of \KPW\ dealt with problems
related to the index. The paper \kpwcite{DouWoj:ALE} with Douglas is a landmark since it
is the starting point of a whole decade seeing a lot of papers focusing
on the $\eta$--invariant and the $\zeta$--determinant.
I was told that it came as an almost unbelievable surprise for the 
mathematical community when $\eta$--function and $\eta$--invariant for 
Dirac operators on compact manifolds {\em with} boundary were established in 
\kpwcite{DouWoj:ALE}, since until 
then the $\eta$-invariant was only established for closed manifolds 
and considered solely as a natural {\em correction} term associated 
to index problems on manifolds with boundary and living exclusively 
on the boundary.

The paper \kpwcite{DouWoj:ALE}
 already contained one of the major analytical
tools which has been refined and exploited ever since: the adiabatic method
(see Section \ref{adiabatic} below). 
\end{sloppy}

There is a variant of the problems mentioned above
which I would like to point out. 
Suppose that $M$ is a closed manifold which is \emph{partitioned}
by a separating hypersurface $\Sigma\subset M$. I. e. there
are compact manifolds with boundary $Y,X$ such that
\footnote{This is a situation which is typical for surgery
theory in which we would have $\Sigma=S^k\times S^l$,
$Y=S^k\times D^l,$ where $S^k$ denotes the unit sphere in $\R^{k+1}$
and $D^k$ denotes the unit disc in $\R^k$.
}
\begin{equation} \label{partmfld}
M=Y\cup_\Sigma X.
\end{equation}

After having chosen appropriate boundary conditions $P^X, P^Y$ for $D$
on $X,Y$ we have three versions of $D$: $D_{P^X}, D_{P^Y}$ and
the essentially self--adjoint operator $D$ on the closed manifold $M$.
In a sense we have "$D=D_{P^Y}\cup D_{P^X}$" and it is natural
to ask how the spectral invariants of $D, D_{P^X},$ and $D_{P^Y}$
are related. \KPW\ and his collaborators have provided us with spectacular
results on this problem.

\subsection{The basic framework of \bvp s for Dirac type operators}

Let us be a bit more specific now and describe the basic set--up of \bvp s
for Dirac type operators as we understand it today. 

Let $X$ and $D$ be as before. We assume that $D$ is an operator of
\emph{Dirac type}. That is in local coordinates
\begin{equation} 
   D^2= - g^{ij}(x)I_{\operatorname{rank} E} \frac{\pl^2}{\pl x_i\pl x_j}+ \text{ lower order
   terms. }
\end{equation}
This is the most general notion of Dirac operator. The leading symbol of 
$D$
\begin{equation}
     \sigma_D(x,df)=i [D,f]_x,\quad f\in\cinf{X},
\end{equation}
induces a \emph{Clifford module} structure on $E$. That is we put for $v\in T_x
X$ \footnote{The Riemannian metric provides us with the ``musical'' isomorphisms
$\flat:T_xM \to T_x^*M$ and $\sharp=\flat^{-1}$.}
\begin{equation} 
     c(v):=\frac 1i \sigma_D(x,v^\flat).
\end{equation}
Then $c(v)^2=-g(v,v)$ and hence by the universal
property of the Clifford algebra, $c$
extends to a section of the bundle $\operatorname{Hom}(Cl(TX,g),\End E)$
of algebra--homomorphisms between the bundle of Clifford--algebras
$Cl(TX,g)$ and the endomorphism bundle $\End E$. This gives $E$ the structure
of a Clifford--module. 

If we choose a Riemannian connection $\nabla$ on $E$ we can form the
Dirac operator $D^\nabla$ on $E$ which is locally given by 
\begin{equation}\label{dnabla}
     D^\nabla=\sum g^{ij} c(\bigl(\frac{\pl}{\pl x_i}\bigr)^\flat) \nabla_{\frac{\pl}{\pl
     x_j}}.
\end{equation}

In the terminology of Boo{\ss}--Bavnbek and Wojciechowski 
\kpwcite{BooWoj:EBP} such operators are called ``generalized
Dirac operators''.
The operators $D^\nabla$ and $D$ obviously have the same leading symbol,
hence 
\begin{equation}
    D=D^\nabla+V
\end{equation}
with $V\in\ginf{X,\End E}$.

Next we have to take the boundary of $X$ into account. We fix a
diffeomorphism from a collar $U$ of the boundary onto
$N:=[0,\eps)\times \Sigma$. Then we may choose a unitary
transformation $\Phi$ from $L^2(U,E)$ onto the product Hilbert space
$L^2([0,\epsilon),L^2(\Sigma,E))$. The operator $\Phi D\Phi^{-1}$ which,
by slight abuse of notation, will again be denoted by $D$ then takes the form
\begin{equation} \label{collarrep}
    D|N=J\bigl(\frac{d}{dx}+B(x)\bigr)+V(x)
\end{equation}
where $J\in\ginf{\Sigma,\End E}$ is a unitary reflection 
($J^2=-I, J^*=-J$), $V\in\cinf{[0,\eps),\ginf{\Sigma,\End E}}$
and $(B(x))_{0\le x\le \eps}$ is a smooth family of first order formally
self--adjoint differential operators on the closed manifold $\Sigma$
(called the tangential operator). 

Replacing $B(x)$ by $B(x)+J^{-1}V(x)$ be obtain alternatively
\begin{equation} \label{collarrep1}
    D|N=J\bigl(\frac{d}{dx}+ \widetilde B(x)\bigr)
\end{equation}
at the expense that now $\widetilde B(x)$ has only self--adjoint leading symbol.

We emphasize that $J$ is independent of $x$
and that \eqref{collarrep} holds for all operators of Dirac type
(Br\"uning and Lesch \kpwcite{BruLes:BVP}, Lemma 1.1).
The representation \eqref{collarrep}
of a generalized Dirac operator is crucial for the geometry of their
\bvp s. In the existing literature, one could sometimes get the impression
that for \eqref{collarrep} to hold one needs that $D$ is the Dirac operator
of a Riemannian connection on $E$ as in \eqref{dnabla} or even a compatible
Dirac operator. 

Furthermore, for many results to be presented below only the following
properties of $D$ will be needed:
\begin{enumerate}
\item $D$ is first order formally self--adjoint elliptic,
\item $D$ has the form \eqref{collarrep} near the boundary,
\item $D$ has the \emph{unique continuation property}.
\end{enumerate}
Properties of Dirac operators which are related to Clifford algebras
will more or less play no role.

$D$ is formally self--adjoint. That is for sections $f,g\in\ginf{X,E}$
we have Green's formula
\begin{equation}\label{green}
    \scalarround{Df}{g}-\scalarround{f}{Dg}=-\int_\Sigma
    \scalar{f}{g}_{E_x}\dvol(x).
\end{equation}
In order to obtain an unbounded self--adjoint operator in $L^2(X,E)$
we have to impose appropriate boundary conditions.

For a pseudo--differential orthogonal projection\footnote{
This is not a big loss of generality.
It can be shown that if the boundary operator has closed range then the
boundary condition may be represented by an orthogonal projection}
\[
P:L^2(\Sigma,E)\to L^2(\Sigma,E)
\]
we define $D_P$ to be the differential
expression $D$ acting on the domain \eqref{domain}.

\begin{definition}\label{wellposed}
\textup{1.} In the notation of \eqref{collarrep} we abbreviate
$B_0:=B(0)$ and denote by $P_+(B_0)$ the orthogonal projection
onto the positive spectral subspace of $B_0$. This is a pseudo--differential
operator of order $0$. Its principal symbol is denoted by $\sigma_{P_+(B_0)}$.

\textup{2.}
The boundary condition defined by $P$ is called \emph{well--posed}
if for each $\xi\in T_x^*\Sigma\setminus\{0\}$ the principal
symbol $\sigma_P(\xi)$ of $P$ maps 
$\range \sigma_{P_+(B_0)}(\xi)$
bijectively onto $\range \sigma_P(\xi)$.
\end{definition}

This is Seeley's definition of well--posedness \kpwcite{See:TPO}. If $P$ is
well--posed then $D_P$ has nice properties.

\begin{proposition} Let $P$ be well--posed. Then $D_P$ is a Fredholm operator
with compact resolvent. Moreover it is regular in the sense
that if a distributional section $u$ of $E$  satisfies
$Du\in L^2_s(X,E)$ and $P(u|\Sigma)=0$ then $u\in L^2_{s+1}(X,E)$,
$s\ge 0$.
\end{proposition}

It turns out that for Dirac type operators this notion of regularity already
characterizes the class of well--posed boundary conditions 
as was shown by Br\"uning and Lesch \kpwcite{BruLes:BVP}.

So far we have basically presented the status of affairs from the point 
of view of classical elliptic theory.

\section{The early work on spectral flow and the general linear conjugation
problem}

\kpwcite{Woj:SFG,BooWoj:DSEI,BooWoj:DSEII}

\bigskip

The early papers \kpwcite{Woj:SFG,BooWoj:DSEI,BooWoj:DSEII} (in part with
Boo{\ss}) on the general
linear conjugation problem are fundamental for our todays understanding
of the structure of \bvp s of Dirac type operators.
The linear conjugation problem is the natural generalization of the classical
Riemann Hilbert problem to elliptic operators (cf. \kpwcite{BooWoj:EBP}, Sec. 26).

Consider a partitioned  manifold $M=Y\cup_\Sigma X$ as in \eqref{partmfld}
and let
\begin{equation} 
    D=\left[ \begin{matrix} 0 & D_-\\
                            D_+& 0\end{matrix}\right]
\end{equation}
be a super-symmetric Dirac operator. That is the bundle $E=E^+\oplus E^-$
is $\Z_2$--graded and $D$ is odd with respect to this grading.

In a collar $N=(-\eps,\eps)\times\Sigma$ of $\Sigma$ we write
$D$ in the form \eqref{collarrep1} and hence we get for $D_+$
\begin{equation} \label{collarrep2}
    D_+=\sigma \bigl(\frac{d}{dx}+B(x)\bigr)
\end{equation}
where $\sigma\in\ginf{\Sigma,\Hom(E^+,E^-)}$ is unitary (and independent of $x$)
and $(B(x))_{-\eps\le x\le \eps}$ is a smooth family of elliptic differential
operators with self--adjoint leading symbol.  

Furthermore, let $\Phi\in\ginf{\Sigma,\Aut(E)}$ be a unitary bundle
automorphism\footnote{\KPW\ originally treated more generally $\Phi$'s
which cover a diffeomorphism of $\Sigma$. Then multiplication by $\Phi$
is a Fourier integral operator.}
of $E$ which is even with respect to the grading.
Multiplication by $\Phi$ is a pseudo--differential operator of order $0$ which
we denote by the same letter. We assume that $\Phi$ commutes with the
leading symbol of $B(x)$. As a consequence the operator $\Phi B \Phi^{-1}-B$ is
of order $0$ and $\Phi P_+(B(x))-P_+(B(x))\Phi$ is of order $-1$ and thus
acts as a compact operator on $L^2(\Sigma,E^+)$.

We introduce a local boundary value problem by letting the differential
expression $D_+$ act on
\begin{equation}\label{lcp}
    \dom(D_+^\Phi):=\bigsetdef{(u_1,u_2)\in L^2_1(Y,E^+)\oplus L^2_1(X,E^+)}{
   u_1|\Sigma=\Phi u_2|\Sigma}.
\end{equation}

From Green's formula \eqref{green} on derives
\begin{equation} 
     (D_+^\Phi)^*=D_-^{\sigma\Phi\sigma^*}
\end{equation}
and thus
\begin{equation}
     D^{\Phi\oplus \sigma\Phi \sigma^*}
          = \left[ \begin{matrix} 0 & D_-^{\sigma\Phi\sigma^*}\\
                            D_+^\Phi& 0\end{matrix}\right]
          =\left[ \begin{matrix} 0 & (D_+^\Phi)^*\\
                            D_+^\Phi& 0\end{matrix}\right].
\end{equation}
One can show that $D^{\Phi\oplus \sigma \Phi\sigma^*}$
is a realization of a local elliptic boundary value problem.
Introducing the \emph{Cauchy data spaces} 
\begin{equation}
     N(D_+,X):=\bigsetdef{u|\Sigma}{u\in L^2_{1/2}(\Sigma,E^+), D_+u=0}
\end{equation}
we find
\begin{equation} \label{indexlcp}
   \begin{split} 
   \ind D_+^\Phi&=\dim\bigl((\Phi N(D_+,X))\cap N(D_-,Y)\bigr)\\
                & \quad -\dim\bigl((J\Phi^*J^* N(D_-,X))\cap N(D_-,Y)\bigr).
   \end{split}
\end{equation}

Before we can state the main result on the linear conjugation problem
we need to elaborate a bit more on the Cauchy data spaces.

\subsection{\Calderon\ projector and the smooth self--adjoint Grassmannian}

\begin{definition}\label{defCald} The (orthogonalized) \Calderon\ projector $C(D,X)$
is the orthogonal projection onto the Cauchy data space $N(D,X)$.
\end{definition}

There is a little subtlety here. The natural construction of the \Calderon\
projector via the invertible double (cf. \kpwcite{BooWoj:EBP}, Sec. 12) gives a
pseudo--differential (in general non--orthogonal)
projection onto the Cauchy data space. It is an orthogonal projection if
$D$ is in product form (cf. \eqref{productform} below) near the boundary. Of course, for any projection 
there is an orthogonal projection with the same image and using the 
results of Seeley \kpwcite{See:CPE} it follows that

\begin{proposition} The orthogonalized \Calderon\ projector $C(D,X)$ is
a pseudo--differential operator of order $0$. Its leading symbol coincides
with the leading symbol $\sigma_{P_+(B_0)}$ of $P_+(B_0)$.
\end{proposition}

The pseudo--differential properties of the \Calderon\ projector
had been developed by \Calderon\ \kpwcite{Cal:BVP}
and Seeley \kpwcite{See:SIB}. 
In \kpwcite{BooLesZhu:IP} we will show that the orthogonalized \Calderon\
projector can be constructed from a natural boundary value problem
on the disconnected double $X\coprod X$.
For brevity we will address the orthogonalized \Calderon\ projector
just as \Calderon\ projector.

The in my view most important observation of the papers 
\kpwcite{BooWoj:DSEI,BooWoj:DSEII} is the fact that the Cauchy
data spaces are Lagrangian. To explain this note that on the Hilbert space
$L^2(\Sigma,E)$ we have the symplectic form
\begin{equation} \label{symplectic}
     \omega(f,g):=-\scalarround{Jf}{g}.
\end{equation}
This claim may be somewhat bewildering since $L^2(\Sigma,E)$ is firstly
a complex vector space and secondly infinite--dimensional. Nevertheless,
$\omega$ is a non--degenerate skew--adjoint sesqui--linear form and it
turns out that it makes perfectly sense to talk about Lagrangians, symplectic
reductions, Maslov indices etc. The only difference is that, due to the
infinite--dimensionality, Fredholm conditions come into play. This is
a fascinating story and an elaboration would definitely need more space.
For some basics cf. Kirk and Lesch \kpwcite{KirLes:EIM}, Sec. 6.
We state explicitly what Lagrangians are in $L^2(X,E)$.

\begin{lemma}\label{lagrangians}
A subspace $L\subset L^2(X,E)$ is Lagrangian if and only if $L^\perp=J(L)$.
\end{lemma}

The following is basically a consequence of Green's formula \eqref{green}.

\begin{proposition}
A realization $D_P$ of a boundary condition
is a symmetric operator if and only if $\range P$ is an isotropic 
subspace of $L^2(X,E)$. Moreover, if $P$ is well--posed then $D_P$
is self--adjoint if and only if $\range P$ is Lagrangian.
\end{proposition}

The following Theorem was proved first in \kpwcite{BooWoj:DSEI}:

\begin{theorem} Let $X$ be a compact Riemannian manifold with boundary
and let $D$ be a Dirac type operator on $X$. Then the Cauchy data
space of $N(D,X)$ is a Lagrangian subspace of $L^2(X,E)$ with respect to
the symplectic structure \eqref{symplectic} induced by Green's form.
\end{theorem}

This theorem is not only beautiful. It is of fundamental importance. We are
now able to describe spaces of well--posed \bvp s as Grassmannian spaces:

\begin{definition} Let $\cP$ be the space of all pseudo--differential orthogonal
projections acting on $L^2(\Sigma,E)$.

The pseudo--differential Grassmannian
$\Gr_1(B_0)$ is the space of $P\in\cP$ such that
\begin{equation}\label{difference}
    P-P_+(B_0) \text{ is of order } -1.
\end{equation}
The space of $P\in\cP$ such that the difference $P-P_+(B_0)$ is smoothing
is denoted by $\Gr_\infty(B_0)$.

Finally the self--adjoint (smooth) pseudo--differential Grassmannian
$\Gr_p^*(B_0)$
is the space of $P\in \Gr_p(B_0), p\in\{1,\infty\}$, whose image is
additionally Lagrangian.
\end{definition}

Since $P_+(B)$ and $C(D,X)$ have the same leading symbol \eqref{difference}
may be replaced by 
\begin{equation}
  P-C(D,X) \text{ is of order } -1.
\end{equation}

Hence $P$ and $C(D,X)$ also have the same leading symbol
and thus it is obvious from the Definition \ref{wellposed} that
the boundary condition given by $P$ is well--posed. 

Furthermore, since the difference of any two elements $P,Q\in\Gr_1(B_0)$
is compact they form a \emph{Fredholm pair}, that is
\begin{equation}
    PQ:\range Q\longrightarrow \range P
\end{equation}
is a Fredholm operator. The index of this Fredholm operator is
denoted by $\ind(P,Q)$. We have
\begin{equation}
     \ind(P,Q)=\dim(\ker P\cap\range Q)-\dim(\range P\cap\ker Q).
\end{equation}

\subsection{The main theorem on the general linear conjugation problem}

We are now in a position to state the main result on the general linear
conjugation problem.

\begin{theorem} The index of the linear conjugation problem \eqref{lcp}
is given by
\[\begin{split}
     \ind D^\Phi&= \ind\bigl(I-C(D_+,Y),\Phi C(D_+,X)\Phi^{-1}\bigr)\\
       &= \ind D+\ind\bigl(C(D_+,X)-\Phi C(D_+,Y)\bigr)\\
       &= \ind D+\ind\bigl(P_+(B_0)-\Phi P_-(B_0)\bigr).
  \end{split}
\]
\end{theorem}

There would be much more to say. This index theorem is related to a lot.
It is a generalization of the classical Riemann Hilbert problem on the 
complex projective line. It is related to the spectral flow and to the
index of generalized Toeplitz operators. 

I will not go into that. But let me say that the papers
\kpwcite{Woj:SFG,BooWoj:DSEI,BooWoj:DSEII} contain much more.
They provide a comprehensive presentation of the spectral flow and
its topological meaning, Fredholm pairs, and the construction of the \Calderon\
projector. Also it is proved that $P_+(B_0)$ is a pseudo--differential
operator.

\section{The $\eta$--invariant}

\kpwcite{DouWoj:ALE,Woj:AEII,Woj:AEIII,LesWoj:IGA,Woj:ZDA}

\bigskip

Let us start with some general remarks on $\eta$-- and $\zeta$--functions.
Let $T$ be an unbounded self--adjoint operator in the Hilbert space $H$. Assume
that $T$ has compact resolvent such that the spectrum of $T$ consists of a 
sequence of eigenvalues 
\[
|\gl_1|\le |\gl_2|\le \ldots \text{ (repeated according to
their finite multiplicity)}
\] 
with $|\gl_n|\to\infty$. If $\gl_n$ satisfies
a growth condition
\begin{equation} 
    |\gl_n|\ge C n^{\ga},
\end{equation}
for some $\ga>0$ then we can form the holomorphic functions
\begin{equation}\begin{split}
     \eta(T;s)&:= \frac{1}{\Gamma(\frac{s+1}{2})}\int_0^\infty
	  t^{(s-1)/2}\tr\bigl(T e^{-tT^2}\bigr)dt\\
        &= \sum_{\gl\in\spec T\setminus\{0\}} |\gl|^{-s} \sign \gl\\
        &= \tr\bigl(T|T|^{-s-1}\bigr),\quad \Re s>\frac 1\ga,
		\end{split}\label{etafctn}
\end{equation}
and
\begin{equation}\begin{split}
     \zeta(T;s)&:= \sum_{\gl\in\spec T\setminus\{0\}} \gl^{-s}\\
        &= \tr\bigl(T^{-s}\bigr),\quad \Re s>\frac 1\ga.
		\end{split}
\end{equation}
If $T$ is non-negative then $\zeta(T;s)$ is also a Mellin transform
similar to the first equality in \eqref{etafctn}
\begin{equation}
     \zeta(T;s)=\frac{1}{\Gamma(s)}\int_0^\infty t^{s-1}\tr \bigl( e^{- t T^2}\bigr)dt.
     \label{zetamellin}
\end{equation}
For general $T$ the function $\zeta(T;s)$ can still be expressed in terms
of Mellin transforms using the formula 
\begin{equation}\label{zetanonselfadjoint} 
     \zeta(T;s)= \frac 12 \bigl(\zeta(T^2;s/2)+\eta(T;s)\bigr)+e^{-i\pi
s}\frac 12\bigl(\zeta(T^2;s/2)-\eta(T;s)\bigr).
\end{equation}

Up to a technical point the existence of a short time asymptotic
expansion of $\tr\bigl(T e^{-tT^2}\bigr), \tr\bigl(e^{-t T^2}\bigr)$ and the meromorphic
continuation of the functions $\zeta(T;s), \eta(T;s)$ is equivalent
(cf. Br\"uning and Lesch \kpwcite{BruLes:EIN}, Lemma 2.2, for the precise statement).

If $T$ is an elliptic operator on a closed manifold then it
follows from the celebrated work of Seeley \kpwcite{See:CPE}
that $\eta(T;s), \zeta(T;s)$ extend meromorphically to $\C$
with a precise description of the location of the poles
and their residues.

If $\eta(T;s)$ is meromorphic at least in a half plane containing
$0$ one defines the
\emph{$\eta$--invariant} of $T$ as 
\begin{equation} \begin{split} 
   \eta(T)&:=\frac{1}{2\pi i} \oint_{|s|=\eps} \frac{\eta(T;s)}{s} ds\\
          &=\text{ constant term in the Laurent expansion at } 0\\
          &=:\eta(T;0).
		 \end{split}
\end{equation}
In many situations one can even show that $\eta(T;s)$ is regular at $0$.
The $\eta$--invariant was introduced in the celebrated work
of Atiyah, Patodi and Singer
\kpwcite{AtiPatSin:SARI,AtiPatSin:SARII,AtiPatSin:SARIII}
as a boundary correction term in an index formula for manifolds
with boundary.

We return to manifolds with boundary and consider again a compact 
Riemannian manifold $X$ with boundary $\pl X=\Sigma$ and
a formally self--adjoint operator of Dirac type acting on
the hermitian vector bundle $E$.

From now on we assume that $D$ is in \emph{product form} near the boundary.
That is in the collar $N=[0,\eps)\times \Sigma$ of the boundary
$D$ takes the form
\begin{equation}\label{productform}
       D|N=J\bigl(\frac{d}{dx}+ B\bigr)
\end{equation}
with $J,B$ as in \eqref{collarrep} and such that
$B$ is independent of $x$. The formal self--adjointness
of $D$ and $B$ then implies
\begin{equation}\label{anticommuting}
     JB+BJ=0.
\end{equation}

The next Theorem guarantees the existence of the $\eta$--invariant
and the $\zeta$--determinant on the smooth self--adjoint Grassmannian.

\begin{theorem}\kpwcite{Woj:ZDA}\label{WojZDA} For $P\in\sGr(B)$ the functions $\eta(D_P;s),\zeta(D_P;s)$
extend meromorphically to a half--plane containing $0$
with poles of order at most $1$.
Furthermore, $0$ is not a pole and $\zeta(D_P;0)$ is independent of
$P$.
\end{theorem}

Let me say a few words about the strategy of proof. As pointed out before
we have to prove  short time asymptotic expansions for $\tr\bigl(D_P e^{-tD_P^2}\bigr)$
and  $\tr\bigl(e^{-tD_P^2}\bigr)$.
Duhamel's principle\footnote{A big word for something very simple:
the method of variation of the constant for first order inhomogeneous
ordinary differential equations.} allows to separate the interior
contributions and the contributions coming from the boundary.
Namely, let $\varphi\in \cinfz{[0,\eps)}$ be a cut--off function
with $\varphi\equiv 1$ near $0$. Extend $\varphi$ by $0$ to a smooth
function on $X$.

Let $\widetilde D$ be any elliptic extension of $D$ to a closed manifold\footnote{
The existence of such a $\widetilde D$ is not essential for the following
result but it simplifies the exposition. For Dirac type operators
we can choose $\widetilde D$ to be the invertible double.}
and let $D_{P,0}$ be the model operator $J\bigl(\frac{d}{dx}+ B\bigr)$
on the cylinder $[0,\infty)\times \Sigma$ with boundary condition $P$
at $\{0\}\times \Sigma$. Then
\begin{equation}\label{tracesplit}\begin{split}
    \tr\bigl(D_Pe^{-t D_P^2}\bigr)&=   \tr\bigl(\varphi D_{P,0}e^{-t D_{P,0}^2}\bigr)+\\
             &\quad  \tr\bigl((1-\varphi)\widetilde D e^{-t \widetilde D^2}\bigr)+O(t^K),\quad
             t\to 0+
		\end{split}
\end{equation}
for any $K>0$. 

By local elliptic analysis the second term in \eqref{tracesplit} has
a short time asymptotic expansion \kpwcite{Gil:ITH2E}, Lemma 1.9.1.
So one is reduced to the treatment of the model operator $D_{P,0}$.
For the Atiyah--Patodi--Singer problem $P=P_+(B)$ there
are explicit formulas for $e^{-t D_{P_+,0}^2}$ from which the asymptotic
expansion can be derived using classical results on special functions.
Finally, for $P\in \sGr(B)$ the operator $D_{P,0}$ can be treated as 
a perturbation of the APS operator $D_{P_+,0}$ \kpwcite{Woj:ZDA}.

A completely different approach by Grubb \kpwcite{Gru:TEP} leads to the
generalization of Theorem \ref{WojZDA} to all well--posed \bvp s.

\subsection{The adiabatic limit}\label{adiabatic}

Let us explain the result of \kpwcite{DouWoj:ALE,Woj:AEII,Woj:AEIII}
on the adiabatic limit of the $\eta$--invariant. We start with
a partitioned manifold $M=Y\cup_\Sigma X$. Then we stretch the neck
by putting

\[\begin{split}\label{adiabaticmfld} 
   X_R&=[0,R]\times \Sigma \cup_{\{R\}\times\Sigma} X,\\
   Y_R&= [-R,0]\times \Sigma\cup_{\{-R\}\times\Sigma} Y,\\
   M_R&=Y_R \cup_{\{0\}\times\Sigma} X_R.
  \end{split}
\]

Denote by $\widetilde \eta (D,M_R)$ the reduced $\eta$--invariant
of $D$ on $M_R$ and by $\widetilde \eta(D_P,X_R)$ the reduced
$\eta$--invariant of $D_P$ on $X_R$. 

\begin{theorem}\label{adiabaticlimiteta}
We have
\begin{equation}\label{ale} \begin{split}
\lim_{R\to\infty} \tilde \eta(D,M_R)&\equiv\lim_{R\to\infty}\tilde
\eta(D_{I-P_+(B)},Y_R)\\
    &\quad+\lim_{R\to\infty}\tilde \eta(D_{P_+(B)},X_R)\mod \Z.
    \end{split}
\end{equation}
\end{theorem}

We should be a bit more specific about the meaning of $P_+(B)$ here.
The positive spectral projection of $B$ is Lagrangian if and only
if $B$ is invertible. If $B$ is not invertible then one has to fix
a Lagrangian subspace of the null space of $B$. So whenever a Lagrangian
is needed we choose $P_+(B)$ such that 
\[
1_{(0,\infty)}(B)\le P_+(B) \le 1_{[0,\infty)}(B). 
\]
That this is possible follows from the Cobordism
Theorem (cf. \kpwcite{DouWoj:ALE} or Lesch and Wojciechowski \kpwcite{LesWoj:IGA}).

In \kpwcite{DouWoj:ALE} it was shown that the $\eta$--invariant makes
sense for generalized Atiyah--Patodi--Singer boundary conditions, i.e.
for $D_{P_+(B)}$. Moreover, it was shown that
$\lim_{R\to\infty}\tilde \eta(D_{P_+(B)},X_R)$ exists. The limit can be
interpreted as the $\eta$--invariant of the operator $D$ on the manifold
with cylindrical ends $X_\infty$. The full strength of Theorem
\ref{adiabaticlimiteta} was proved in \kpwcite{Woj:AEII,Woj:AEIII}.
In fact the ($\textup{mod } \Z$ reductions) of the ingredients of formula
\eqref{ale} do not depend on $\Z$ as was observed by W. M\"uller
\kpwcite{Mul:EIM}. In this way we obtain the gluing formula for the
$\eta$--invariant for the boundary condition $P_+(B)$. The following
generalization to all $P\in \sGr(B)$ is worked out in \kpwcite{Woj:ZDA}.


\begin{theorem} Let $M=Y\cup_\Sigma X$ be a partitioned manifold
and let $D$ be a Dirac type operator which is in product form
in a collar of $\Sigma$. Then for $P\in\sGr(B)$
\begin{equation}
     \tilde\eta(D,M)\equiv \tilde\eta(D_P,X)+\tilde\eta(D_{I-P},Y)\mod\Z.
\end{equation}
\end{theorem}


There is even a formula if $I-P$ is replaced by a general $Q\in\sGr(-B)$.
This is an extension of a formula for the variation of the $\eta$--invariant
under a change of boundary condition from \kpwcite{LesWoj:IGA}, cf. also
Theorem \ref{ML-S4.2} below.

\bigskip
Because of its importance let us look briefly at the method of proof.

The first observation is that the heat kernel of the model operator
$D=J\bigl(\frac{d}{dx}+ B\bigr)$ on the cylinder $\R\times \Sigma$ is explicitly
known since $D^2$ is just a direct sum of one--dimensional Laplacians
$-\frac{d^2}{dx^2}+b^2$. Let $\cE_{\cyl}(t;x,y)$ be this cylinder heat
kernel. Furthermore, denote by $\cE_R(t;x,y)$ the heat kernel of $D$ on
the stretched manifold $M_R$.

Next one chooses $R$--dependent cut--off 
functions $\phi_{j,R},\psi_{j,R}, j=1,2,$ as follows:
\begin{align*}
\psi_{2,R}(x)&=
\begin{cases}
0 \quad &\text{if $|x|\le 3R/7,$}\\
1 &\text{if $|x|\ge 4R/7,$}
\end{cases}\\  
\psi_{1,R}&=1-\psi_{2,R}.
\end{align*}
Finally, choose $\phi_{j,R}$ such that $\phi_{j,R}\psi_{j,R}=\psi_{j,R}$.
Then paste the heat kernel $\cE_R$ on $M_R$ and the cylinder heat kernel
to obtain the kernel
\begin{equation}
   Q_R(t;x,y)=\phi_{1,R}(x)\cE_{\cyl}(t;x,y)\psi_{1,R}(y)
             +\phi_{2,R}(x)\cE_R(t;x,y)\psi_{2,R}(y).
\end{equation}
Then Duhamel's principle yields
\begin{equation}
  \cE_R(t)=Q_R(t)+\cE_R\# C_R(t),
\end{equation}
where $\#$ is a convolution and $C_R$ is an error term.

It seems that not much is gained yet. The point is that
Douglas and Wojciechowski \kpwcite{DouWoj:ALE} could show
that in the adiabatic limit the error term is negligible
in the following sense:


\begin{theorem} There are estimates
\[\begin{split}
     \|\cE_R(t;x,y)\|&\le c_1 t^{-\dim X/2}e^{c_2 t} e^{-c_3d^2(x,y)/t},\\
      \|(\cE_R\# C_R)(t;x,x)\|&\le c_1 e^{c_2 t} e^{-c_3R^2/t}
  \end{split}
\]
with $c_1,c_2,c_3$ independent of $R$.
\end{theorem}


Note that this result is much more than e.g. \eqref{tracesplit}.
For the $\eta$-- and $\zeta$--determinant the full heat semigroup
contributes. It is astonishing that nevertheless in the adiabatic
limit the full integrals from $0$ to $\infty$ in \eqref{eta} and
\eqref{zeta} split into contributions coming from the cylinder
and from the interior of the manifold.

\section{The relative $\eta$--invariant and the relative $\zeta$--determinant}

\kpwcite{LesWoj:IGA,ScoWoj:ZDQ}

\bigskip

Recall from Theorem \ref{WojZDA} that for $P\in\sGr(B)$ the
$\zeta$--function $\zeta(D_P;s)$ is regular at $0$. 
One puts
\begin{equation}
     \detz D_P:=\begin{cases}\exp\bigr(-\zeta'(D_P;0)\bigr), &
0\not\in\spec D_P,\\
                            0, & 0\in\spec D_P.
               \end{cases}
\label{ML-G4.1}
\end{equation}
In view of \eqref{zetanonselfadjoint} 
and Theorem \ref{WojZDA} a straightforward
calculation shows for $D_P$ invertible
\begin{equation}
    \detz D_P= \exp\Bigl( i\frac{\pi}{2} \bigl(\zeta(D_P^2;0)-\eta(D_P)\bigr)
                  -\frac 12 \zeta'(D_P^2;0)\Bigr).
    \label{ML-G4.2}
\end{equation}
We emphasize that the regularity of $\eta(D_P;s)$ and $\zeta(D_P;s)$ at
$s=0$ is essential for \eqref{ML-G4.2} to hold.
\eqref{ML-G4.2} shows that the $\eta$--invariant
is related to the phase of the $\zeta$--determinant and that
in general 
\[
(\detz D)^2\not= \detz (D^2).
\]

The natural question which arises at this point is


\begin{problem}\label{problem} How does $\detz(D_P)$ depend on $P\in\sGr(B)$?
\end{problem}


The answer to this problem has a long history. 
Since the only joint paper of Wojciechowski and myself
deals with an aspect of the problem I take the liberty
to add a few personal comments. In 1992 I was a Postdoc
at University Augsburg.
At that time the paper \kpwcite{DouWoj:ALE}
had just appeared and the gluing formula for the $\eta$--invariant
was in the air. Still much of our todays understanding
of spectral invariants for Dirac type operators on manifolds with
boundary was still in its infancy. When Gilkey visited
he posed  
a special case of the Problem \ref{problem}. If the tangential operator is
not invertible there is no canonical Atiyah--Patodi--Singer
boundary condition for $D$. The positive spectral projection
of $B$ is not in $\Gr_\infty^*(B)$. Rather one has to
choose a Lagrangian subspace $V\subset \ker B$ and put
\[
P_V:=1_{(0,\infty)}(B)+\Pi_V, 
\]
where $\Pi_V$ denotes the orthogonal
projection onto $V$. Then $P_V\in\Gr_\infty^*(B)$.
The boundary condition given by $P_V$ is called a generalized
Atiyah--Patodi--Singer boundary condition. Gilkey asked
how the eta--invariant depends on $V$.

I did some explicit calculations on a cylinder which let
me guess the correct formula. However, I did not know how
to prove it in general. Somewhat later Gilkey sent me
a little note of Krzysztof dealing with the same problem.
He urged us to work together. I was just a young postdoc
and I felt honored that Krzysztof, whose papers I already
admired, quickly agreed. Except writing papers with my
supervisor this was my first mathematical collaboration.
It was done completely by fax and email; Krzysztof and I
met for the first time more than a year after the paper
had been finished.


In \kpwcite{LesWoj:IGA} Krzystof and I proved a special
case of the following result. The result as stated is
a consequence of the Scott--Wojciechowski Theorem as was shown
in \kpwcite{KirLes:EIM}, Sec. 4. The Scott--Wojciechowski
Theorem will be explained below.


\begin{theorem}\label{ML-S4.2}  Let $P,Q\in \sGr(B)$. Then
\begin{equation} \etab(D_P)-\etab(D_{Q})
\equiv \log\detf(\Phi(P)\Phi(Q)^*)\mlmod \Z.\label{ML-G4.5}
\end{equation}
\end{theorem}


If $P$ or $Q$ is the \Calderon\ projector then \eqref{ML-G4.5}
is even an equality \kpwcite{KirLes:EIM}.

\medskip

The general answer to Problem \ref{problem}
given by Scott and Wojciechowski \kpwcite{ScoWoj:ZDQ} is just beautiful.
To explain their result we need another bit of notation.
Recall that $J$ defines the symplectic form on $L^2(\Sigma,E)$
\eqref{symplectic}. Let 
\[
E=E_i\oplus E_{-i}
\]
be the decomposition of $E$
into the eigenbundles of $J$. If $P\in\sGr(B)$ then 
\[
L=\range P \subset L^2(\Sigma,E)
\]
is Lagrangian and from Lemma \ref{lagrangians}
one easily infers that the restrictions of the orthogonal projections
$\Pi_{\pm i}=\frac{1}{2i}(i\pm J)$ onto $E_{\pm i}$ map $L$
bijectively onto $L^2(\Sigma,E_{\pm i})$ and 
\[
\Phi(P):= \Pi_{-i}\circ (\Pi_i|E_i)^{-1} 
\]
is a unitary operator from
$L^2(\Sigma,E_i)$ onto $L^2(\Sigma,E_{-i})$. For $P$ we then have the
formula
\begin{equation}
    P=\frac 12\begin{pmatrix} I & \Phi(P)^*\\ \Phi(P) &I\end{pmatrix}.
\end{equation}
For $P,Q\in\sGr(B)$ the operator $\Phi(P)^*\Phi(Q)-I$ is smoothing and
hence $\Phi(P)^*\Phi(Q)$ is of determinant class.


With these preparations, the Scott--Wojciechowski theorem
reads as follows.


\begin{theorem}\label{ML-S4.1} Let $P\in\sGr(B)$ and let $C(D,X)$ be
the orthogonalized \Calderon\ projector. Then
\begin{equation}    \detz(D_P)=\detz(D_{C(D,X)})
       \detf\bigl(\frac{I+\Phi(C(D,X))\Phi(P)^*}{2}\bigr).
       \label{ML-G4.3}
\end{equation}
\end{theorem}

\section{Adiabatic decomposition of the $\zeta$--determinant}

\kpwcite{KliWoj:ACT,ParWoj:ADZ,ParWoj:STA}

\bigskip
When the gluing formula for the $\eta$--invariant had been established
it was Krzysztof's optimism that eventually lead to a similar result
for the $\zeta$--determinant. The author has to admit that he was
an unbeliever: I could not see why a reasonable analytic surgery formula
for the $\zeta$--determinant should exist. Well, I was wrong.
A fruitful collaboration of J. Park and \KPW\ eventually proved
that the adiabatic method, which originally had been developed in the 
paper \kpwcite{DouWoj:ALE}, was even strong enough to prove an adiabatic
surgery formula for the $\zeta$--determinant.

Consider again the adiabatic setting $M_R, X_R, Y_R$ as in 
\eqref{adiabaticmfld}. In order not to blow up the exposition too much
I will not present the result in its most general form. Rather I will make
the following technical assumptions:

\begin{enumerate}
\item The tangential operator $B$ is invertible.
\item The $L^2$--kernel of $D$ on $X\cup[0,\infty)\times \Sigma$ and
$Y\cup[0,\infty)\times \Sigma$ vanishes.
\end{enumerate}

Then the adiabatic surgery theorems for the Laplacians read as follows:


\begin{theorem} Let $\Delta_{\pm,R,d}$ be the Dirichlet extension of $D^2$ on
  $X_R, Y_R$ resp.; $D_R$ denotes the operator $D$ on $X_R$. Then
\[
\lim_{R\to\infty} \frac{\detz D_R^2}{\detz \Delta_{+,R,d} \detz\Delta_{-,R,d}}
= \sqrt{\detz B^2}.\]
\end{theorem}


\begin{theorem} Let 
$D_{+,R,P_+}$, $D_{-,R,P_-}$ be the operator $D$ with Atiyah--Patodi--Singer
boundary conditions on $X_R, Y_R$ resp. Then
\[  
\lim_{R\to\infty} \frac{\detz D_R^2}{\detz D^2_{+,R,\Pi_{>}}{\detz D^2_{-,R,\Pi_{<}}}}
= 2^{-\zeta'(B^2,0)}.\]
\end{theorem}

These technical assumptions mentioned above
were removed in Park and Wojciechowski \kpwcite{ParWoj:STA}. For details the reader should consult
loc. cit.

\smallskip

Finally, the ``adiabatic'' results on the zeta--determinants obtained by
Park and Wojciechowski are not adiabatic any more. Loya and Park 
\kpwcite{LoyPar:GPS}
showed that most of those results (and more) are true without stretching.
Krzysztof P. Wojciechowski did have different (and charming) 
ideas how to remove
stretching of the cylinders. Unfortunately, his serious illness
did not allow him to fill all the details and finish the paper. 

\newcommand{\Toappear}{to appear in}

\begin{thebibliography}{10}

\bibitem{AtiPatSin:SARI}
M.~F. Atiyah, V.~K. Patodi, and I.~M. Singer, \emph{{S}pectral asymmetry and
  {R}iemannian geometry {I}}, Math. Proc. Cambridge Philos. Soc. \textbf{77}
  (1975), 43--69.

\bibitem{AtiPatSin:SARII}
M.~F. Atiyah, V.~K. Patodi, and I.~M. Singer, \emph{{S}pectral asymmetry and
  {R}iemannian geometry {II}}, Math. Proc. Cambridge Philos. Soc. \textbf{78}
  (1975), 405--432.

\bibitem{AtiPatSin:SARIII}
\bysame, \emph{Spectral asymmetry and {R}iemannian geometry {III}}, Math. Proc.
  Cambridge Philos. Soc. \textbf{79} (1976), 71--99.

\bibitem{BooWoj:DSEI}
B.~Boo{\ss} and K.~P. Wojciechowski, \emph{Desuspension and splitting
  elliptic symbols {I}}, Ann. Global Anal. Geom. \textbf{3} (1985), 337--383.

\bibitem{BooWoj:DSEII}
\bysame, \emph{Desuspension and splitting elliptic symbols {II}}, Ann. Global
  Anal. Geom. \textbf{4} (1986), 349--400.


\bibitem{BooLesZhu:IP}
B.~Boo{\ss}-Bavnbek, M.~Lesch, and C.~Zhu, \emph{In preparation}.

\bibitem{BooWoj:EBP}
B.~Boo{\ss}-Bavnbek and K.~P. Wojciechowski, 
\emph{Elliptic boundary problems for {D}irac operators}, Birkh\"auser,
Basel, 1993.

\bibitem{BruLes:EIN}
J.~Br{\"u}ning and M.~Lesch, \emph{On the eta--invariant of certain non--local
  boundary value problems}, Duke Math. J. \textbf{96} (1999), 425--468,
  dg-ga/9609001.

\bibitem{BruLes:BVP}
\bysame, \emph{On boundary value problems for {D}irac type operators: {I}.
  {R}egularity and self--adjointness}, J. Funct. Anal. \textbf{185} (2001),
  1--62, math.FA/9905181.

\bibitem{Cal:BVP}
A.~Calder{\'o}n, \emph{Boundary value problems for elliptic equations},
  Outlines of the joint Soviet--American symposium on partial differential
  equations (Novosibirsk), 1963, pp.~303--304.

\bibitem{DouWoj:ALE}
R.~G. Douglas and K.~P. Wojciechowski, \emph{Adiabatic limits of the
  $\eta$-invariants the odd--dimensional {A}tiyah--{P}atodi--{S}inger problem},
  Comm. Math. Phys. \textbf{142} (1991), 139--168.

\bibitem{Gil:ITH2E}
P.~Gilkey, \emph{Invariance theory, the heat equation, and the
  {A}tiyah--{S}inger index theorem}, 2. ed., CRC Press, Boca Raton, 1995.

\bibitem{Gru:TEP}
G.~Grubb, \emph{Trace expansions for pseudodifferential boundary problems for
  {D}irac--type operators and more general systems}, Ark. Mat. \textbf{37}
  (1999), 45--86.

\bibitem{KirLes:EIM}
P.~Kirk and M.~Lesch, \emph{The $\eta$--invariant, {M}aslov index, and spectral
  flow for {D}irac--type operators on manifolds with boundary}, Forum Math.
  \textbf{16} (2004), 553--629, math.DG/0012123.

\bibitem{KliWoj:ACT}
S.~Klimek and K.~P. Wojciechowski, \emph{Adiabatic cobordism theorems for
  analytic torsion and $\eta$--invariant}, J. Funct. Anal. \textbf{136} (1996),
  269--293.

\bibitem{LesWoj:IGA}
M.~Lesch and K.~P. Wojciechowski, \emph{On the $\eta$--invariant of generalized
  {A}tiyah--{P}atodi--{S}inger boundary value problems}, Illinois J. Math.
  \textbf{40} (1996), 30--46.

\bibitem{LoyPar:GPS}
P.~Loya and J.~Park, \emph{On the gluing problem for the spectral invariants of
  dirac operators}, To appear.

\bibitem{Mul:EIM}
W.~M{\"u}ller, \emph{Eta invariants and manifolds with boundary}, J.
  Differential Geom. \textbf{40} (1994), 311--377.

\bibitem{ParWoj:ADZ}
J.~Park and K.~P. Wojciechowski, \emph{Adiabatic decomposition of the
  $\zeta$--determinant of the {D}irac {L}aplacian {I}. {T}he case of invertible
  tangential operator}, Comm. Partial Differential Equations \textbf{27}
  (2002), 1407--1435, (with an Appendix by Y. Lee).

\bibitem{ParWoj:STA}
\bysame, \emph{Scattering theory and adiabatic decomposition of the
  $\zeta$--determinant of the {D}irac {L}aplacian}, Math. Res. Lett. \textbf{9}
  (2002), 17--25.

\bibitem{ScoWoj:ZDQ}
S.~G. Scott and K.~P. Wojciechowski, \emph{The $\zeta$--determinant and
  {Q}uillen determinant for a {D}irac operator on a manifold with boundary},
  Geom. Funct. Anal. \textbf{10} (1999), 1202--1236.

\bibitem{See:TPO}
R.~T. Seeley, \emph{Topics in pseudo--differential operators}, C.I.M.E.,
  Conference on pseudo--differential operators 1968, Edizioni Cremonese, Roma,
  1969, pp. 169--305.

\bibitem{See:SIB}
\bysame, \emph{Singular integrals and boundary value problems}, Amer. J. Math.
  \textbf{88} (1966), 781--809.

\bibitem{See:CPE}
\bysame, \emph{Complex powers of an elliptic operator}, Proc. Sympos. Pure
  Math. \textbf{10} (1967), 288--307.

\bibitem{Woj:SFG}
K.~P. Wojciechowski, \emph{Spectral flow and the general linear conjugation
  problem}, Simon Stevin \textbf{59} (1985), 59--91.

\bibitem{Woj:AEII}
\bysame, \emph{The additivity of the $\eta$--invariant. {T}he case of an
  invertible tangential operator}, Houston J. Math. \textbf{20} (1994),
  603--621.

\bibitem{Woj:AEIII}
\bysame, \emph{The additivity of the $\eta$--invariant. {T}he case of a
  singular tangential operator}, Comm. Math. Phys. \textbf{109} (1995),
  315--327.

\bibitem{Woj:ZDA}
\bysame, \emph{The $\zeta$-determinant and the additivity of the
  $\eta$-invariant on the smooth, self-adjoint {G}rassmannian}, Comm. Math.
  Phys. \textbf{201} (1999), 423--444.

\end{thebibliography}


\providecommand{\bysame}{\leavevmode\hbox to3em{\hrulefill}\thinspace}
\providecommand{\MR}{\relax\ifhmode\unskip\space\fi MR }
\providecommand{\MRhref}[2]{%
  \href{http://www.ams.org/mathscinet-getitem?mr=#1}{#2}
}
\providecommand{\href}[2]{#2}

\editorsnote
\end{document}